\documentclass{amsart}
\font\smc=cmcsc10 %%%%%%%%%%
\font\ssmc=cmcsc8 %%%%%%%%%%
\newcommand\sep{,\ }  %%%%
\newtheorem{lemma}{Lemma}
\newtheorem{theorem}{Theorem}
\newtheorem*{example}{Example}
\newtheorem*{corollary}{Corollary}

\newcommand\NN{\mathbb{N}}
\newcommand\QQ{\mathbb{Q}}
\newcommand\RR{\mathbb{R}}
\newcommand\ZZ{\mathbb{Z}}
\newcommand\ZV{\ZZ V}
\newcommand\AG{\overline{\NN\mathcal{A}}_G}
\newcommand\Fb{\overline{F}}
\newcommand\fp{\mathfrak{p}}
\newcommand\fq{\mathfrak{q}}

\begin{document}
\title{Toric ideals of normalized graph algebras} % integral closures from graphs}
\author{Peter M.\ Johnson}
\address{Departamento de Matem\'atica, Universidade Federal da Bahia, 
Salvador, Brazil.}
%%%%%%%%%%%%%%%%%  ?? How does elsarticle handle current address ??
\curraddr{\ssmc Depto de Matem\'atica, UFPE, %% Univ.\ Federal de Pernambuco, 
50670-901 Recife--PE, Brazil.}
\email{peterj@dmat.ufpe.br}

\begin{abstract}
A graph-theoretic method, simpler than existing ones, is used to
characterize the minimal set of monomial generators for the
integral closure of any algebra of polynomials generated by quadratic
monomials.  The toric ideal of relations between these generators
is generated by a set of binomials, defined graphically.  The
spectra of the original algebra and of its integral closure turn out
to be canonically homeomorphic.
\end{abstract}

\maketitle

% \begin{keyword} %%%%%%%%
\noindent{\smc Keywords:}  
{Polynomial ring\sep integral closure\sep quadratic monomial\sep toric ideal.}
{\smc 2010 MSC Class}: Primary: 13B22; Secondary: 13F20, 05C99.  
%% 2010 same as 2000
%% Does elsarticle use Primary/Secondary ??  

\section{Introduction: monomial algebras from graphs}

We hope to draw the attention of non-specialists to a relatively accessible area of
commutative algebra with interesting questions that can be explored using a
rich variety of techniques, notably combinatorial ones.
Our focus is on a single problem, the structure of the integral closures (also known
as normalizations) of certain algebras of polynomials defined from graphs.
Questions about objects constructed from graphically-defined algebras can be regarded as 
questions about graphs, but little is known about when one can expect answers naturally 
expressible in graph-theoretical language.  To understand integral closures of the algebras, 
it is desirable to have not only a good set of generators but also detailed information 
about the relations they satisfy, something that has until now been lacking.
As we show, inherently algorithmic methods from graph theory quickly determine the
generators, simplifying known proofs, and lead to an especially clear description of when 
polynomials in generators are equal.  A few graphically-defined kinds of relations between 
generators will be found to suffice.
% generate all others.

Arguments based on graphical configurations, of the kind in
Section~8.7 of \cite{Vr} or in related articles such as \cite{SVV}, are our principal tool.
We draw more heavily than usual on such methods, and diagrams should be drawn as needed.
Although the sequence of ideas is new, the methods are so elementary
that parts of arguments inevitably recur in the literature.
For example, it is useful to know how even closed walks in graphs can be split.
This problem also arose, for a different reason, in \cite{BGS}, a recent article
in the same general area.  A more efficient approach to the solution is presented below.

Some background from commutative algebra will be briefly described.
All rings considered are commutative algebras over a field $k$.  The study of integral 
closures of powers and symbolic powers of ideals of polynomial rings is of interest, 
but (using Rees algebras) it usually suffices to concentrate only on subalgebras.
Villarreal's monograph \cite{Vr} describes a rich and still developing theory, 
involving for example polyhedral combinatorics and hypergraphs, that arises  
from the study of objects related to monomials. This subarea of commutative algebra has 
close connections with the subarea of algebraic geometry concerned with toric varieties.

For any subalgebra $F$ of a ring $R$ of polynomials (or Laurent polynomials) over a field $k$,
the integral closure $\Fb$ of $F$ can be regarded as a subalgebra of $R$.
If $F$ is generated by (Laurent) monomials then $\Fb$, as a $k$-vector space,
has a basis consisting of those monomials $f$ in the radical of $F$ 
that are quotients $f = g/h$ of monomials in $F$.
This basic result, which opens the door to combinatorial methods, is widely used,
and can be proved using an elementary idea from polyhedral combinatorics, 
but we could not find a single article or book that identifies its source.
Of five experts consulted, only Prof.\ Hibi supplied the reference: it is an immediate 
consequence of an auxiliary result (Prop.~1) in a classic article by Hochster [Ho] from 1972.
Briefly, Hochster characterizes the monomially-generated algebras that are normal 
(integrally closed), and one sees from his conditions that the normalization $\Fb$ of $F$ 
must coincide with the ring obtained by adjoining to $F$ the monomials in $\Fb$, 
which are just those described above.

Specializing further, to algebras $F$ generated by quadratic monomials $x_ix_j$ or $x_i^2$,
permits more detailed results, such as those of \cite{SV}, which are stated in terms of
graphs $G$ on $\{x_1, \dots, x_n\}$ whose edges correspond to the given generators of $F$.
In this special case, the problem of giving a graphical description for monomial generators 
of the integral closure $\Fb$, using the basis described above as a starting point, is 
trivially equivalent to a purely graph-theoretic problem that will soon be stated.
Explicit generators were first found by Simis, Vasconcelos, Villarreal  \cite{SVV} and by
Ohsugi, Hibi \cite{OH}, in articles submitted only two months apart, which cover much the 
same ground despite minor differences in assumptions and forms of results.
The approach of Ohsugi and Hibi emphasized the theory of integral polytopes in spaces
$\RR^d$, and their set of generators can be seen to be minimal, while
Simis, Vasconcelos, and Villarreal used little more than elementary polynomial algebra,
except in their Proposition~2.1, and chose to allow some redundant generators.
Section~8.7 of \cite{Vr} contains a partly simplified exposition, and that book
contains a variety of related topics.

Our more elementary approach to finding generators provides additional
information about the monomials in $\Fb$.
Such ideas are exploited to describe binomial generators for the ideal of polynomial relations
between the given generators of $\Fb$ (its toric ideal).  A description has until now been
lacking, the main difficulty being that, while the problem is of purely algebraic interest,
its solution seems to require an unusually intensive use of methods from discrete mathematics.
It is too much to expect a natural definition of an irredundant set of generating 
relations, but a relatively small generating class will be described. 
We hope that our solution and methods will help stimulate more refined investigations of
rings defined from graphs.

\section{Generators for integral closures via weighted graphs}

Throughout $G = (V,E)$ is a graph, allowing loops but not multiple edges, with vertex set 
$V = V_G = \{x_1, \dots, x_n\}$.  
Traditionally, the graphs used here are finite, although this restriction is inessential.  
Some fairly standard graph-theoretic terminology 
is presupposed.  Much use is made of walks $w$ in graphs, recorded as ordered lists 
of adjacent vertices or of the corresponding edges.  Their size is measured by 
the number of edges, and the parity (odd or even) of walks will often be important.
Walks are either open or closed.
For walks that are closed (the last vertex equals the first), the vertices will be regarded 
as being in a cyclic rather than a linear order, so that the walk is unchanged under the 
obvious rotations.  Paths are walks in which all vertices except possibly the first 
and last are different.  Circuits are closed paths.

In this area, rather than studying the rings $F$ and $\Fb$ directly, it is more common to
work with the associated semigroups of monomials, written additively.
Let $\ZV$ (embedded in $\QQ V$) denote the free abelian group generated by 
the vertices $V$ of $G$.
Subsets $X$ of $V$ will be identified with their characteristic functions in $\ZV$.
An ordinary edge $\{x_i,x_j\}$ is represented as $x_i+x_j$, while $2x_i$ is used for a loop 
at $x_i$.  Under addition, the subset $\mathcal{A}_G$ of images of the edges of $G$ generates 
a monoid $\NN\mathcal{A}_G$, an object often studied.  
By a convention going back to \cite{Ho}, this monoid is usually called a semigroup.
Elements of $\ZV$ can sometimes 
arise by assigning non-integral weights to edges in subgraphs, where the weight at a vertex $v$ 
is the sum of the weights of the edges at $v$, with loops at $v$ counted twice.  For example, 
every circuit with edges of weight $\frac{1}{2}$ induces weight 1 on each of its vertices.
This weighted circuit, if even, can also be written as a sum (with all weights 1) of certain edges.

The above problem about the integral closure $\Fb$ of a ring $F$ defined from a graph $G$, translated 
into additive language, becomes that of finding generators for the semigroup called $\AG$ formed by
the elements of $\ZV$ that are linear combinations of edges with nonnegative rational coefficients.  
Curiously, the problem of describing relations between the generators has remained untouched until 
now, although a good set of generators for $\AG$ has long been known, from \cite{OH} or \cite{SVV}.
We first describe these generators, then give a short new proof that they suffice.

Following \cite{OH}, except that $G$ is allowed to be disconnected, an {\sl exceptional pair}
in $G$ is an induced subgraph $H$ consisting of two vertex-disjoint odd circuits, perhaps loops,
that lie in the same connected component of $G$.  
These coincide with H-configurations, as defined in \cite{SV}, when $G$ is connected.
It seems best to work with such objects rather than with similar but less restrictive ones known
as bow ties, where $H$ need not be induced.  If $H$ has an odd circuit with a chord in $G$ or 
has vertex-disjoint circuits that can be joined by an edge in $G$, it is well-known and easy 
to verify that $H$ is a sum in $\ZV$ of edges plus (except in the last case) some $H'$,
where $H'$ is an exceptional pair supported on a proper subset of the vertices on 
circuits of $H$.  As explained after Def.~2.5 of \cite{SV},  
$H$ (in $\ZV$) is a weighted sum, with coefficients $\pm 1$, of edges.
Also $2H$, but not $H$ itself, is a positive integral linear combination of edges, 
so $H \in \AG$.

\smallskip

With this terminology and notation, the result on generators for $\Fb$ becomes:

\begin{theorem}  The edges and exceptional pairs of the graph $G$ form
the unique minimal generating set of $(\AG,+)$.
\end{theorem}

\begin{proof}
Each $f \in \AG$ is a positive rational linear combination of the edges in some set $E_f$.
By an easy result known as Carath\'eodory's Theorem, it can be assumed that these edges are
linearly independent in $\QQ V$.  Indeed, if $E_f$ is dependent, one can replace
$f$ repeatedly by $f$ plus a multiple of a nontrivial linear relation, making some
coefficients of vectors in $E_f$ zero and keeping the others positive.
We work with pairs $(f,E_f)$, with $E_f$ independent, and show how to generate $f$ by
inducting on $|E_f|$, which is zero only when $f=0$.
Note that the subgraph $G_f = (V, E_f)$ of $G$ cannot have any even circuits,
so its connected components are trees or contain a unique odd circuit.  
Also, $f$ determines uniquely the weights on the edges in $E_f$.

Separation of these edge weights into integer and fractional parts decomposes $f$
within $\AG$ as $f = g+h$, where $g$ is a positive integral sum of edges, 
$E_h \subseteq E_f$, and $h$ is a sum, with weights strictly between 0 and 1, 
of the edges in $E_h$.  This reduces the argument to the case that $f=h$.
Also, from the definition of $\AG$, all vertex weights $f(v)$ are integers.
Thus, in $G_f$, no vertex lies on a unique edge, so this graph
must be a disjoint union of odd circuits of $G$ (allowing loops) on whose
vertices $f$ has weight 1.

By assumption, $f$ is also a $\ZZ$-linear combination of edges. It follows that, 
for each connected component $K$ of $G$, the sum of the weights of $f$ at vertices 
in $K$ is even, so $K$ must contain an even number of circuits of $G_f$.  These 
circuits can then be grouped into pairs each of which forms, in $G$, an exceptional 
pair or a non-induced subgraph of the kind handled in an earlier discussion.

Thus $\AG$ is generated (as a monoid) by edges and exceptional pairs.
None of these can be rewritten as a sum of other generators of $\AG$,
so such elements form the unique minimal generating set.
\end{proof}

\section{Minimal relations and moves between monomials}

The aim of this section is to prepare the ground for a later description of the 
polynomial relations satisfied by edges and exceptional pairs in the integral 
closure $\Fb$, where as always $F$ is defined from some graph $G$ on the vertices
$x_1, \dots, x_n$.
Abstractly, one uses two new kinds of variables, here simply called edges and pairs
where confusion is unlikely.
These map to certain monomials (products of the vertices involved) in
$\Fb$, regarded as a subring of $R = k[x_1, \dots, x_n]$.
The polynomial relations sought are of the form $p = q$, where $p$ and $q$ are
monic monomials, usually called \emph{words}, in edges and pairs.  These in turn 
give binomials $p-q$ which clearly suffice to generate the toric ideal of $\Fb$.
The relations of equality can also be thought of as defining a congruence relation
on the commutative semigroup of words in the new variables.  We wish to describe 
a few kinds of graphically-defined relations that generate this congruence.

It is convenient to introduce auxiliary variables, here called \emph{cycles},
associated with induced odd circuits in the graph $G$, with weight $\frac{1}{2}$ 
on the edges and thus weight 1 on the vertices of the circuit. 
This unconventional use of `cycle' should be noted.
Each pair variable converts into a product of two cycle variables.
It is much simpler to analyze words in edges and cycles, then use the results  
to solve the original problem of describing relations between words in edges and pairs.
In arguments, cycles (variables) often stand in for odd circuits (subgraphs),
 
Passing from $\AG$ to a larger subsemigroup $\AG^+$ of $\ZV$ generated by edges 
and cycles, each word gives an edge-weighted subgraph of $G$ which in turn gives 
an element of $\AG^+$.
We intend to describe the relations of equality between words in edges and cycles
that are minimal, in the sense that no proper nonempty subsum yields a valid equality.
All relations of equality are clearly sums of minimal ones, often in several ways.
Cancellation of common terms, which would in general reduce greatly the number
of generating relations needed, is an operation not permitted here.
The notion of minimality has one weakness: it does not exploit the transitivity 
of equality, so the complete list of minimal relations may be highly redundant. 
The following trivial example illustrates the problem.

\begin{example}
Let $G$ be formed from three non-adjacent loops by adding edges connecting these to a 
new vertex. The element of $\AG$ having all vertex weights 1 can be expressed in three
ways as a sum of an edge and a pair. This gives three minimal equations in generators of 
$\AG$, each a consequence of the other two.  Thus none of these relations is irredundant.
\end{example}

We now analyze a purely graph-theoretic concept relevant to the minimality of certain 
relations. Whenever an even closed walk (in cyclic order) does not split at some vertex
into two consecutive even closed walks, but is of the form $w_1w_2w_3w_4$,
where $w_1$ and $w_3$ are walks from a vertex $x_1$ to $x_2$, possibly with $x_1 = x_2$,
then the closed walks $w_1w_2$ and $w_2w_3$ are odd, so $w_1$ followed by the reverse 
of $w_3$ is an even closed walk. Thus the walk splits in a less direct way, and
the equation from alternate edge sums of the original closed walk is
a sum of equations from two other such walks.

All other even closed walks that do not split, but have a repeated vertex $v$, are
of the form $w = w_1w_2$, where $w_1$ and $w_2$ begin and end with $v$.
Then $w_1$ and $w_2$ are odd closed walks with only $v$ in common.  
One can remove $v$ (fusing the edges at $v$ in pairs) to create two vertex-disjoint 
even closed walks from which $w$ can be reconstructed by gluing the smaller walks 
together at a new vertex placed within an edge of each walk.
In this way, every even closed walk that does not split can be constructed recursively 
from even circuits, here including the usually forbidden cases of circuits with 
no vertices or with one edge, repeated twice.  An equivalent result on splitting, 
via a different approach, is presented in Section~3 of \cite{BGS}.

The following three kinds of equality relations between words, valid in $\AG^+$, 
will turn out to include all minimal ones.
It is often preferable to think in terms of directed moves that alter parts of words in 
generators, or the related edge-weighted graph, without affecting the image in the semigroup.

\medskip %%%
% \begin{description} 
% \item[Rotation moves.]
{\bf Rotation moves. \ }
The weighted sum of edges in an even closed walk, with weights alternating between 
1 and -1, vanishes. This gives a relation of equality between two edge sums with the
same number of terms, counting multiplicities if edges repeat.
It is not hard to see that the relation is minimal precisely when the walk
does not split in the way discussed immediately above. The corresponding rotation move 
starts with edge-weights 0 and 1 that alternate along the walk, and adjusts them by $\pm 1$, 
thus interchanging the edge-weights 0 and 1 without altering the induced vertex-weights.

It is harmless to introduce a new variable for each even closed circuit, representing 
the sum with weights $\frac{1}{2}$ of the edges in the circuit.
This permits new moves, here called \emph{half-rotations}, to convert between a 
circuit variable and either of the two sums of alternate edges in the circuit.
In effect, each edge weight is adjusted by $\pm \frac{1}{2}$.
Then each rotation can be realized by two half-rotations. 

\medskip %%%
{\bf Cycle-destroying moves. \ }
For even closed walks formed from the edges in an exceptional pair and in a walk $w$ 
between the two circuits in the pair, traversed twice in the larger even closed walk,
the sum of alternate (either odd or even) edges in $w$ plus the two cycles from the pair 
can be rewritten as a sum of some of the edges in the larger even closed walk.  
An associated move (there are two) can be thought of as a new kind of half-rotation 
without a new variable.  
As the name suggests, one usually applies a cycle-destroying move in the direction 
that eliminates its pair of cycles.

A single cycle is never a sum of edges with integral weights,
as is clear from the parities of sums of vertex weights in the odd circuit.
Thus a cycle-destroying relation is minimal 
precisely when the walk $w$ between the circuits contains no even closed subwalk. 
In the more degenerate case of two odd circuits too close to form an exceptional
pair, maybe even equal, the cycle sum becomes an edge sum.  As an extreme case, each loop
has a double role as an edge $e$ and a cycle $c$, where $c+c = e$.

\medskip %%%
{\bf Cycle-shifting moves. \ }
The associated relations are those valid in $\AG^+$ 
that have, on each side of the equality, a single cycle 
variable plus a sum of edges, allowing repetitions.
In addition, to ensure minimality, the sides must have no terms in common, 
not even if one side is first altered by a rotation move.

As will be seen in the proof of Theorem~2, all cycle-shifts arise from configurations 
in $G$ with the following properties, where cycles (really the underlying circuits) 
and edges may overlap or repeat.  Given two different cycles and a fixed set of 
open walks in $G$, the extreme vertices of all these walks are different, and are
the vertices that lie on one cycle but not on both.
Also, each of these walks has odd (resp.\ even) length when it begins and ends in 
the same cycle (resp.\ different cycles).  
It follows that two types, called positive and negative signs, can be assigned 
to the generators (edges or cycles) so that signs alternate along the walks, 
and the two cycles have different signs.   The sums of generators of each type 
are then equal, forming a minimal relation.

\medskip %%%
Some redundancy involving the first two kinds of moves can easily be eliminated.

\begin{lemma}
Every rotation move can be realized by a sequence of cycle-destroying moves and 
rotation moves on non-split even closed circuits that are not the union of two smaller 
even closed circuits minus a common edge.
\end{lemma}

\begin{proof}
Inducting on the number of edges in the relevant walks,
one can suppose that the even closed walk $w$ associated with
the given rotation move does not split. 
If $w$ is a circuit that violates the last condition of the lemma,
the rotation can be carried out via two rotation moves
on smaller circuits.  For $w$ not a circuit, the earlier analysis 
of non-split walks shows that $w$ must contain
at least two odd circuits $C_1$ and $C_2$, say $w = w_1C_1w_2C_2$,
where $w_1w_2$ forms a smaller even closed walk, possibly null.
Now consider a configuration with weights $\frac{1}{2}$ on the edges of $C_1$ and $C_2$, 
and weights $1$ and $-1$ alternating along the walk $w_1w_2$.  
One now chooses a cycle-destroying move $M_1$ (resp.\ $M_2$) using $w_1$ (resp.\ $w_2$) to 
connect the two circuits.  With the right choice of these moves, the desired rotation move 
can now be obtained by applying $M_1$ in reverse followed by $M_2$.
\end{proof}

Alternatively, cycle-shifts and certain kinds of half-rotations generate all moves.

\section{Relations between monomials in standard form}

Two words in edges and cycles will of course map to the same element of $\AG^+$
precisely when their formal difference induces the weight 0 on all vertices of $G$, 
but what is really wanted is a proof that sequences of known moves can transform one
word into the other.  Our method will provide an efficient algorithm for doing this.
Ideas involving flux in weighted graphs guide all the arguments.

A word in edges and cycles is said to be in \emph{standard form} if 
the corresponding edge-weighted graph admits no cycle-destroying moves.  
The next theorem is a striking illustration 
of the value of converting words into standard form.
Here words are treated as sets, counting multiplicities.

\begin{theorem}
Suppose $w_1$ and $w_2$ are disjoint nontrivial words in standard form that represent the 
same element of $\AG^+$.  Then some subword of $w_1$ can be moved, by a single rotation or 
cycle-shift, to a subword of $w_2$.  Every repeated selection of such subwords on the unused 
parts of $w_1$ and $w_2$ terminates with a partition of $w_1$ and one of $w_2$, with moves 
that transform the parts of $w_1$ into the parts of $w_2$.
\end{theorem}

\begin{proof}
After reducing as much as possible, it can be assumed that no moves on subwords
can be made between $w_1$ and $w_2$ and there are no common generators, but the
words are not null.
Generators in $w_1$ (resp.\ $w_2$) will be called positive (resp.\ negative),
and one can work with edge weights induced by the formal difference of these words.
All edges in $w_1$ and in $w_2$ (counting positive multiplicities) can be grouped
into walks where consecutive edges alternate in sign, forming each walk successively
to be as long as possible at both ends.  None of these walks can be even and closed, 
as that would allow a rotation move between subwords of $w_1$ and $w_2$.  From vertex 
weights, one can now see that the walks are open, and their extreme vertices (all different) 
are the vertices that lie on exactly one of the cycles of $w_1$ or $w_2$.  The alternating 
walks can then be regarded as starting and ending with an appropriate cycle.  
These are the only walks considered below.  As the words admit no cycle-destroying moves, 
the odd walks begin and end at the same cycle, while the even walks join cycles of 
different signs.  Thus each cycle is joined to an odd number of cycles of the other sign.
Remarkably, this correspondence between cycles will be shown to be a bijection.

To prove this, fix a cycle $C$ of vertices and edges in $G$ corresponding to a cycle variable 
in $w_1$ (say), and for each cycle $C_i$ from $w_2$, let $V_i$ denote the set of vertices
of $C$ at which walks from $C_i$ reach $C$.  Such sets are mutually disjoint.
For distinct cycles $C_i$, $C_j$ from $w_2$, and vertices $x_i \in V_i$, $x_j \in V_j$,
no pair $\{x_i,x_j\}$ can be an edge of $C$, as this would allow a cycle-destroying move
in $w_2$.  Similarly, if $\{x_i,u\}$ and $\{v, x_j\}$ are edges of $C$,
no walk has $u$ and $v$ as its extreme vertices.  More generally, one can start
with a cyclic order (usually not induced by $G$) on some subset of vertices of $C$,
and at each stage until no longer possible remove from the list some pair of vertices $u,v$ 
joined by a walk and create a new cyclic order on the new set, as follows.
If the previous cyclic order can be written as walks $w_1w_2$ (usually not in $G$) from $u$ 
to $v$ to $u$, the new cyclic order consists of $w_1$ followed by the reverse of $w_2$, 
now omitting $u$ and $v$ from these lists of vertices.
By induction, vertices adjacent via such cyclic orders can never
lie in sets $V_i$, $V_j$ that connect to different cycles $C_i$, $C_j$.
At the final step, it becomes clear that only one cycle $C_i$ from $w_2$
is joined to $C$ by walks of the above form.
By symmetry, the correspondence between cycles in $w_1$ and those in $w_2$ is a bijection.
But this yields a cycle-shifting move between $w_1$ and $w_2$, contrary to assumption.

More constructively, the analysis produces a partition of the generators of $w_1$
into supports for rotations and cycle-shifts that transform $w_1$ into $w_2$.
\end{proof}

The following useful conclusions can immediately be drawn.

\begin{corollary}
Words in standard form that represent the same element contain equal numbers of cycle variables.
Such words remain in standard form when rotation or cycle-shifting moves are applied.
\end{corollary}

\section{Consequences for integral closures of graph algebras}

The algebras $\Fb$ are generated by edges and pairs, rather than edges and cycles.
The relations are closely related to the ones already studied, and the details are
routine, so we shall be concise, now using multiplicative language.
To simplify statements and the definition of pairs, $G$ will be assumed to be connected.  
One can, however, easily generalize results by treating each connected component of $G$ separately. 
For pairs, we use variables $p_{ij}$ indexed by unordered pairs of indices that enumerate 
the induced odd circuits in $G$, with the understanding that $p_{ij}$ either represents an 
exceptional pair formed from cycles $C_i$ and $C_j$, or is merely an abbreviation 
for some product of edges obtained from destroying this pair of cycles. 

Several earlier definitions need to be adapted to the new variables.
Rotation moves (or relations) involve only edges, so no change is necessary. 
A cycle-destroying move, which destroys two cycles in an exceptional pair,
translates into a pair-destroying move, and a word not admitting such a move 
is said to be in standard form.  Thus, apart from the change of variables, the concept
of standard form does not change under the translation.  For pair-shifting equations,
it suffices to admit only those obtainable by multiplying a cycle-shifting equation
on both sides by the same cycle variable, then replacing on each side the pair of cycles
$C_i$, $C_j$ by a pair variable $p_{ij}$.
We call such a move a pair-shifting move with a shared cycle.  
 
\begin{theorem}
The toric ideal of the normalized algebra $\Fb$ defined from a graph,
with edges and pairs as generators,
is generated by the relations $p_{ij}p_{kl} = p_{ik}p_{jl}$ and those from 
rotation moves, pair-destroying moves, and pair-shifting moves with a shared cycle.  
\end{theorem}

\begin{proof}[Sketch of proof:]
Words in the new generators induce words in edges and cycles, and the procedure in the previous 
theorem for transforming words into standard form and using moves between words with equal images 
translates immediately to the new situation.  Note that the way that cycles are paired can be 
changed whenever convenient, and two cycle-shifts applied in parallel to a word in standard form 
are now realizable via two pair-shifts, each with a different shared cycle.
\end{proof}

In the present situation, with words formed from edges and pairs, the notion of minimality
is no longer useful.  One way to obtain complicated examples of minimal equations in edges
and pairs is to choose a suitable pairing of cycles in a move (as in Theorem~2)   %%
composed of many disjoint cycle-shifts. 
 
The above analysis of relations may form a starting point for answering more sophisticated
questions about algebraic objects related in some way to $\Fb$.
  
Going in an apparently unexplored direction, we end with elementary observations
about prime spectra of graphically-defined algebras, ignoring the sheaf structure.  
Basic results related to the going-up theorem, found in standard sources 
such as Ch.~5 of Atiyah and MacDonald \cite{AM}, say in spectral language that the 
continuous map $Spec(\Fb) \to Spec(F)$ induced by the inclusion $F \to \Fb$ is a 
closed surjection between (quasi)compact spaces, with discrete fibers.
Thus this is a homeo\-morphism if it is injective.

\begin{theorem}
For $F$ as above, $Spec(\Fb)$ is canonically homeomorphic to $Spec(F)$.
\end{theorem}

\begin{proof}
We fix a prime ideal $\fp$ of $F$, and let $M$ be a sufficiently large field
containing the domain $F/\fp$, so that the one or more prime ideals $\fq$ of $\Fb$
that lie above $\fp$ are the kernels of the homomorphisms $\Fb \to M$ that extend
the quotient map $F \to F/\fp$.    To obtain $M$, it suffices to start with a field $K$
generated by $F/\fp$ and make simultaneous quadratic extensions 
so that $M$ is of the form $K[t_1,\dots,t_m]$ with $[M:K] = 2^m$, where the $t_i$ 
correspond to certain pairs $p_i$ whose squares lie in $F \setminus \fp$.  
The above homomorphisms $\Fb \to M$ are determined by the image of each $p_i$ 
($t_i$ or $-t_i$).   One sees without even needing to invoke Galois theory that, 
over $K$, $M$ has $2^m$ automorphisms (1 if char$(K)=2$).
Then all these homomorphisms $\Fb \to M$ have the same kernel, so there is a
unique prime ideal $\fq$ above $\fp$.
\end{proof}

We mention only in passing that these affine schemes can be glued together to form 
less trivial examples of schemes that do not ramify under normalization.  With this 
in mind, it seems worth describing how the graph $G$ can be used to clarify the topological 
structure of $Spec(F)$.  Using notations such as $F = F_K$ to emphasize the subgraph $K$
of $G$ used, 
let $L_K$ be the algebra of Laurent polynomials over $k$ in the edges of $K$.  A free generating 
set of edges and their inverses is obtained from a spanning forest for $K$ to which, in each 
non-bipartite connected component of $K$, an extra edge is added to form an odd circuit.
With such generators, rings of the form $L_K$ have well-understood prime spectra.

Next, call a subgraph $K$ of $G = (V,E)$ \emph{admissible} if, 
for each even closed circuit in $G$, 
with edges considered to be of two types that alternate around the circuit, 
$K$ contains all edges of one type whenever it contains all of the other type.
Each $\fp \in Spec(F)$ defines an admissible graph whose edges are
$\{ e \in E \mid e \notin \fp\}$.
Conversely, for each admissible subgraph $K$, let $\fp_K$ denote the ideal of $F$ 
generated by the edges not in $K$. Known results about relations in $F$, like our 
more general Theorem~2, give $F_G/\fp_K \cong F_K$.  Thus $\fp_K \in Spec(F)$,
and still gives a prime ideal after adjoining inverses of the edges of $K$.
These conclusions, couched in spectral language, yield:
  
\begin{theorem}
After canonical identifications via localizations, 
the sets $Spec(L_K)$, $K$ an admissible subgraph of $G$,
partition $Spec(F)$ into locally closed subsets.
The points of $Spec(L_K)$ are the prime ideals of $F$ giving the graph $K$, 
and their least member (unique closed point) is $\fp_K$.
\end{theorem}
  
\medskip %%%

%%%%%%%%%%   How can I make elsarticle print the heading REFERENCES ??

\end{document}